\documentclass{article}
\usepackage{amsmath}
\usepackage{amssymb}
\usepackage{amsfonts}

\newtheorem{theorem}{Theorem}

\newtheorem{lemma}[theorem]{Lemma}

\newenvironment{proof}[1][Proof]{\noindent \textbf{#1.} }{\  \rule{0.5em}{0.5em}}

\begin{document}

\author{$^{a}$Zafer \c{S}iar and $^{b}$Refik Keskin \and $\ ^{a}$Bing\"{o}l
University, Department of Mathematics, Bing\"{o}l/TURKEY \and $\ ^{b}$%
Sakarya University, Department of Mathematics, SakaryaTURKEY \and $^{a}$%
zsiar@bingol.edu.tr, $^{b}$rkeskin@sakarya.edu.tr, }
\title{The solutions of the Diophantine equation $(2^{k}-1)(3^{l}-1)=7^{m}-1$%
}
\maketitle

\begin{abstract}
In this paper, we will show that the positive integer triple $(k,l,m)$
satisfying the equation in the title is $(2,1,1).$
\end{abstract}



\section{Introduction}

The Fibonacci sequence $(F_{n})$ is defined by $F_{0}=0,~F_{1}=1,$ and $%
F_{n}=F_{n-1}+F_{n-2}$ for $n\geq 2.$ The Lucas sequence $(L_{n})$, which is
similar to the Fibonacci sequence, is defined by the same recursive pattern
with initial conditions $L_{0}=2,~L_{1}=1$. The terms of the Fibonacci and
Lucas sequences are called Fibonacci and Lucas numbers, respectively.

Let $a,b,c$ be positive integers such that $1<a<b<c.$ Let consider the
Diophantine equation%
\begin{equation}
(a^{k}-1)(b^{l}-1)=c^{m}-1  \label{1.1}
\end{equation}%
in positive integers $k,l,m.$ For $(a,b,c)=(F_{3},F_{4},F_{5}),$ the authors
\cite{Luca-S} proved that the triple $(k,l,m)$ satisfying the equation (\ref%
{1.1})\emph{\ }is only $(2,2,2).$ Later, this problem was extended by Irmak
\cite{Nurettin} to the case $(a,b,c)=(F_{n},F_{n+1},F_{n+2}).$ Irmak showed
that the positive integer quadruples $(n,k,l,m)$ satisfying the equation (%
\ref{1.1})\emph{\ }are $(3,2,2,2)$ and $(5,2,1,2).$

In this paper, we deal with this problem\emph{\ }for Lucas numbers. Namely,
we will prove the following theorem.

\begin{theorem}
\label{T1} The only positive integer triple $(k,l,m)$ satisfying the
equation
\begin{equation}
(2^{k}-1)(3^{l}-1)=7^{m}-1  \label{1.2}
\end{equation}%
is $(2,1,1).$
\end{theorem}

\section{Preliminaries}

Let $\eta $ be an algebraic number of degree $d$ with minimal polynomial
\begin{equation*}
a_{0}x^{d}+a_{1}x^{d-1}+\cdots +a_{d}=a_{0}\prod \limits_{i=1}^{d}\left(
x-\eta ^{(i)}\right) \in \mathbb{Z}[x],
\end{equation*}%
where the $a_{i}$'s are relatively prime integers with $a_{0}>0$ and the $%
\eta ^{(i)}$'s are conjugates of $\eta .$ Then
\begin{equation}
h(\eta )=\frac{1}{d}\left( \log a_{0}+\sum \limits_{i=1}^{d}\log \left( \max
\left \{ |\eta ^{(i)}|,1\right \} \right) \right)  \label{2.1}
\end{equation}%
is called the logarithmic height of $\eta .$ Particularly, if $\eta =a/b$ is
a rational number with $\gcd (a,b)=1$ and $b\geq 1,$ then $h(\eta )=\log
\left( \max \left \{ |a|,b\right \} \right) .$

The following theorem is deduced from Corollary 2.3 of Matveev \cite{Mtv}
and provides a large upper bound for the unknowns $k,l,$ and $m$ in equation
(\ref{1.2}) (also see Theorem 9.4 in \cite{Bgud}).

\begin{theorem}
\label{T2}Assume that $\gamma _{1},\gamma _{2},...,\gamma _{t}$ are positive
real algebraic numbers in a real algebraic number field $\mathbb{K}$ of
degree $D$, $b_{1},b_{2},...,b_{t}$ are rational integers, and
\begin{equation*}
\Lambda :=\gamma _{1}^{b_{1}}\cdots \gamma _{t}^{b_{t}}-1
\end{equation*}%
is not zero. Then
\begin{equation*}
|\Lambda |>\exp \left( -1.4\cdot 30^{t+3}\cdot t^{4.5}\cdot D^{2}(1+\log
D)(1+\log B)A_{1}A_{2}\cdots A_{t}\right) ,
\end{equation*}%
where
\begin{equation*}
B\geq \max \left\{ |b_{1}|,...,|b_{t}|\right\} ,
\end{equation*}%
and $A_{i}\geq \max \left\{ Dh(\gamma _{i}),|\log \gamma _{i}|,0.16\right\} $
for all $i=1,...,t.$
\end{theorem}

During our calculations, we get upper bounds on our variables which are too
large, so we have to reduce them. To do this, we use some results from the
theory of continued fractions. Regarding this, in \cite{duj}, Dujella and
Peth\H{o} gave a result, which is a version of the reduction method based on
the Baker-Davenport Lemma \cite{Baker}. Then, in \cite{bravo3}, the authors
proved the following lemma, which is an immediate variation of the result
due to Dujella and Peth\H{o} from \cite{duj}.

\begin{lemma}
\label{L2}Let $M$ be a positive integer, let $p/q$ be a convergent of the
continued fraction of the irrational number $\gamma $ such that $q>6M,$ and
let $A,B,\mu $ be some real numbers with $A>0$ and $B>1.$ Let $\epsilon
:=||\mu q||-M||\gamma q||,$ where $||\cdot ||$ denotes the distance from $x$
to the nearest integer. If $\epsilon >0,$ then there exists no solution to
the inequality
\begin{equation*}
0<|u\gamma -v+\mu |<AB^{-w},
\end{equation*}%
in positive integers $u,v,$ and $w$ with
\begin{equation*}
u\leq M\text{ and }w\geq \frac{\log (Aq/\epsilon )}{\log B}.
\end{equation*}
\end{lemma}

The above lemma allows us find a lower bound for linear forms in two integer
variables. But, in some cases, we need to find a lower bound for linear
forms in three or more integer variables in the form of

\begin{equation*}
x_{1}\tau _{1}+x_{2}\tau _{2}+\cdots +x_{t}\tau _{t}\text{ with\ }t\geq 3,
\end{equation*}%
where $\tau _{1},\tau _{2},\ldots ,\tau _{t}$ are real numbers and $%
x_{1},x_{2},\ldots ,x_{t}$ are integers. In this case, we use the LLL
algorithm that we describe below.

Put $X=\max \left\{ \left\vert x_{1}\right\vert ,\left\vert x_{2}\right\vert
,\ldots ,\left\vert x_{t}\right\vert \right\} $, $C>(tX)^{t}$ and let $%
\Omega $ be a lattice generated by
\begin{equation*}
b_{j}=e_{j}+\left\lfloor C\tau _{j}\right\rceil e_{t}\text{ for }1\leq j\leq
t-1\text{ and }b_{t}=\left\lfloor C\tau _{t}\right\rceil e_{t},
\end{equation*}%
where $C$ is sufficiently large positive constant.

\begin{lemma}
\label{L4}\emph{\cite[Lemma 5]{dda}} Let $X_{1},X_{2},\ldots ,X_{t}$ be
positive integers such that $X:=\max \left\{ X_{i}\right\} $ and $C>(tX)^{t}$
is a fixed sufficiently large constant. With the above notation on the
lattice $\Omega $, we consider a reduced base $\{b_{i}\}$ to $\Omega $ and
its associated Gram-Schmidt orthogonalization base $\{b_{i}^{\ast }\}$. We
set%
\begin{equation*}
c_{1}:=\max_{1\leq i\leq t}\frac{\left\Vert b_{1}\right\Vert }{\left\Vert
b_{i}^{\ast }\right\Vert },~\theta :=\frac{\left\Vert b_{1}\right\Vert }{%
c_{1}},~S:=\sum\limits_{i=1}^{t-1}X_{i}^{2}\text{ and }T:=\frac{1}{2}\left(
1+\sum\limits_{i=1}^{t}X_{i}\right) .
\end{equation*}%
If the integer $x_{i}$ are such that $\left\vert x_{i}\right\vert \leq X_{i}$
for $1\leq i\leq t$ and $\theta ^{2}>S+T^{2},$ then we have
\begin{equation*}
\left\vert \sum\limits_{i=1}^{t}x_{i}\tau _{i}\right\vert \geq \frac{\sqrt{%
\theta ^{2}-S}-T}{C}.
\end{equation*}
\end{lemma}

The following lemma is useful in the proof of the main theorem.

\begin{lemma}
\label{L3}\emph{\cite[Lemma 2.2]{weger}} Let $a,u\in
%TCIMACRO{\U{211d} }%
%BeginExpansion
\mathbb{R}
%EndExpansion
$ and $0<$ $a<1$. If $\left\vert u\right\vert <a,$ then
\begin{equation*}
\left\vert \log (1+u)\right\vert <\frac{-\log (1-a)}{a}\cdot \left\vert
u\right\vert
\end{equation*}%
and
\begin{equation*}
\left\vert u\right\vert <\frac{a}{1-e^{-a}}\cdot \left\vert
e^{u}-1\right\vert .
\end{equation*}
\end{lemma}

\section{The proof of Theorem \protect\ref{T1}}

Assume that the equation (\ref{1.2}) holds with $k,l,m\geq 1.$ It is clear
that $k$ can not equal to $1.$ Assume that $k\geq 2.$ Then we see that $%
(-1)((-1)^{l}-1)\equiv ((-1)^{m}-1)(mod4)$, and from this, we get%
\begin{equation*}
2\equiv (-1)^{m}+(-1)^{l}(mod4),
\end{equation*}%
which implies that $m$ and $l$ have the same parity. On the other hand,
since $2^{k}-1\equiv 0,1,3(mod7)$ and $(2^{k}-1)(3^{l}-1)=7^{m}-1\equiv
-1(mod7),$ it follows that $3^{l}\equiv 3(mod7),$ which implies that $l=6q+1$
for some integer $q.$ That is, $l$ is odd.

Now, let us record what we have proved.

\begin{lemma}
\label{L5} If the equation (\ref{1.2}) holds with $k,l,m\geq 1,$ then $l$
and $m$ are odd.
\end{lemma}

\begin{lemma}
\label{L6}Assume that the equation (\ref{1.2}) holds. Then $k=2$ or $l=1$ if
and only if $m=1.$
\end{lemma}

%TCIMACRO{\TeXButton{Proof}{\proof} }%
%BeginExpansion
\proof
%EndExpansion
Let $k=2.$ Then, we have $\left( 3^{\frac{l+1}{2}}\right) ^{2}-2=7^{m}$ has
only the solution $m=1$ by Lemma 2.5 given in \emph{\cite{ting}.} Let $l=1.$
Then we have the equation $2^{k+1}-7^{m}=1.$ This equation has no solutions
for $m>1$ by Lemma 2.2 given in \emph{\cite{ting}}. So, $m=1$. Conversely,
if $m=1,$ the proof follows. This completes the proof of Lemma \ref{L6}.

So, from now on, we suppose that $k\geq 3$ and $l\geq 2.$ Thus, $m\geq 3.$
Now, let rearrange the equation (\ref{1.2}) as follows:%
\begin{equation}
7^{m}=2^{k}3^{l}-2^{k}-3^{l}+2.  \label{3.0}
\end{equation}%
So, we get $7^{m}<2^{k}3^{l}.$ Also, since $7^{m}>7^{m}-1\geq 7\cdot
(3^{l}-1)>6\cdot 3^{l}$ and $7^{m}>7^{m}-1\geq 8\cdot (2^{k}-1)>6\cdot 2^{k},
$ it follows that $7^{m}>6\max \left\{ 2^{k},3^{l}\right\} \geq 3\cdot
\left( 2^{k}+3^{l}\right) .$ Thus, we have%
\begin{equation}
3\cdot \left( 2^{k}+3^{l}\right) \leq 6\max \left\{ 2^{k},3^{l}\right\}
<7^{m}<2^{k}3^{l}.  \label{3.1}
\end{equation}%
Moreover, with the help of Mathematica program, it can be seen that
\begin{equation}
2^{k-0.2}<2^{k}-1  \label{3.2}
\end{equation}%
for $k\geq 3$ and
\begin{equation}
3^{l-0.2}<3^{l}-1  \label{3.3}
\end{equation}%
for $l\geq 2.$ Using (\ref{3.1}), (\ref{3.2}) and (\ref{3.3}), we
get\noindent
\begin{equation*}
2^{k-0.2}\cdot 3^{l-0.2}<(2^{k}-1)(3^{l}-1)=7^{m}-1<7^{m}<2^{k}3^{l}.
\end{equation*}%
Taking logarithm in the above inequalities and making necessary
calculations\noindent , we obtain%
\begin{equation}
0.35k+0.56l-1.82<m<0.36k+0.57l<\max \left\{ k,l\right\} .  \label{3.4}
\end{equation}%
\noindent Now, let arrange equation (\ref{1.2}) as follows:\noindent
\begin{equation*}
\frac{2^{k}\cdot 3^{l}}{7^{m}}-1=\frac{3^{l}+2^{k}-2}{7^{m}}.
\end{equation*}%
\noindent If we use (\ref{3.0}) and (\ref{3.1}), we get%
\begin{eqnarray}
0 &<&\frac{2^{k}\cdot 3^{l}}{7^{m}}-1<\frac{2^{k}+3^{l}}{7^{m}}\leq \frac{%
2\max \left\{ 2^{k},3^{l}\right\} }{7^{m}}  \label{3.5} \\
&=&\frac{2\max \left\{ 2^{k},3^{l}\right\} \min \left\{ 2^{k},3^{l}\right\}
}{7^{m}\cdot \min \left\{ 2^{k},3^{l}\right\} }=\frac{2\cdot 2^{k}\cdot 3^{l}%
}{7^{m}\cdot \min \left\{ 2^{k},3^{l}\right\} }  \notag \\
&<&\frac{2\cdot \left( 7^{m}+2^{k}+3^{l}\right) }{7^{m}\cdot \min \left\{
2^{k},3^{l}\right\} }<\frac{2\cdot \left( 7^{m}+7^{m}/3\right) }{7^{m}\cdot
\min \left\{ 2^{k},3^{l}\right\} }<\frac{8/3}{2^{\min \left\{ k,l\right\} }}.
\notag
\end{eqnarray}%
Put $\Lambda _{1}:=1-{{\gamma }_{1}}^{b_{1}}{{\gamma }_{2}}^{b_{2}}{{\gamma }%
_{3}}^{b_{3}},$ where\noindent
\begin{equation*}
{\gamma }_{1}:=2,{\gamma }_{2}:=3,{\gamma }%
_{3}:=7,b_{1}:=k,b_{2}:=l,b_{3}:=-m.
\end{equation*}%
If $\Lambda _{1}=0,$ then we get $2^{k}\cdot 3^{l}=7^{m}$, which is
impossible since $k\geq 3$, $l\geq 2$, and $m\geq 3$ are integers. Therefore
$\Lambda _{1}\neq 0.$ Thus, we can apply Theorem \ref{T2} to inequality (\ref%
{3.5}). Note that the numbers ${\gamma }_{1},{\gamma }_{2}$ and ${\gamma }%
_{3}$ are positive real numbers and elements of the field of $\mathbb{%
%TCIMACRO{\U{211a} }%
%BeginExpansion
\mathbb{Q}
%EndExpansion
}.$ The degree of this field is $1.$ So $D=1.$ Moreover, \noindent since
\begin{equation*}
h\left( {\gamma }_{1}\right) =h\left( 2\right) ={\log 2,}
\end{equation*}%
\begin{equation*}
h\left( {\gamma }_{2}\right) =h\left( 3\right) ={\log 3\ ,}
\end{equation*}%
\begin{equation*}
h\left( {\gamma }_{3}\right) =h\left( 7\right) ={\log 7\ ,}
\end{equation*}%
by (\ref{2.1}), we can take\noindent\ $A_{1}:={\log 2}$, $A_{2}:={\log 3}$, $%
{A_{3}:=\log 7}$ in Theorem \ref{T1}. Also, using (\ref{3.4}), we see
that\noindent
\begin{equation*}
\max \left\{ k,l,m\right\} =\max \left\{ k,l\right\} .
\end{equation*}%
Say $B=\max \left\{ k,l\right\} .$ Thus, by Theorem \ref{T1}, we
obtain\noindent
\begin{equation*}
\frac{8/3}{2^{\min \left\{ k,l\right\} }}>\left\vert \Lambda _{1}\right\vert
>exp\left( -1.4{\cdot 30}^{6}{\cdot }3^{4.5}\cdot \left( 1+\log B\right)
\cdot \log {2\cdot \log 3{\cdot \log 7{\ }}\ }\right) ,
\end{equation*}%
\noindent which implies that%
\begin{equation}
\min \left\{ k,l\right\} <3.1\cdot {10}^{11}\cdot \left( 1+\log B\right) .
\label{3.6}
\end{equation}%
Now, we distinguish the proof into two cases.\noindent

\textbf{Case 1: }Let $l\leq k.$ That is $B=k.$ In this case, we have
\begin{equation}
l<3.1\cdot {10}^{11}\cdot \left( 1+\log k\right)   \label{3.7}
\end{equation}%
from (\ref{3.6}). On the other hand, using (\ref{3.0}) and (\ref{3.1}), we
can arrange equation (\ref{1.2}) as follows\noindent :\
\begin{equation}
0<\frac{2^{k}\cdot \left( 3^{l}-1\right) }{7^{m}}-1=\frac{3^{l}-2}{7^{m}}<%
\frac{3^{l}}{7^{m}}<\frac{4/3}{2^{k}}.  \label{3.8}
\end{equation}%
Put\noindent\ $\Lambda _{2}:=1-{{\gamma }_{1}}^{b_{1}}{{\gamma }_{2}}^{b_{2}}%
{{\gamma }_{3}}^{b_{3}}$, where\noindent
\begin{equation*}
{\gamma }_{1}:=2,{\gamma }_{2}:=3^{l}-1,{\gamma }%
_{3}:=7,b_{1}:=k,b_{2}:=1,b_{3}:=-m.
\end{equation*}%
If $\Lambda _{2}=0,$ then\noindent\ we get%
\begin{equation*}
2^{k}\cdot \left( 3^{l}-1\right) =7^{m},
\end{equation*}%
which implies that $k=0.$ But, this is not possible since $k\geq 1.$
Therefore $\Lambda _{2}\neq 0.$ A similar arguments to the above shows that
we can take $D=1,$ $B:=k$, $A_{1}:={\log 2}$, $A_{2}:=l\cdot \log 3>{\log
\left( 3^{l}-1\right) \ }$ and $A_{3}:={\log 7}$ in Theorem \ref{T2}. Thus,
applying Theorem \ref{T2} to inequality (\ref{3.8}), we get\noindent
\begin{equation}
k<3.1\cdot 10^{11}\cdot \left( 1+\log k\right) {\cdot l}.  \label{3.9}
\end{equation}%
\noindent If we substitute the bound of $l$ in (\ref{3.7}) into (\ref{3.9}),
then we get%
\begin{equation*}
k<9.61\cdot 10^{22}\cdot \left( 1+\log k\right) ^{2}.
\end{equation*}%
which implies that\noindent\
\begin{equation}
l\leq k<3.72\cdot 10^{26}.  \label{3.10}
\end{equation}%
\textbf{Case 2: }Let $k<l.$ Then, we have
\begin{equation}
k<3.1\cdot {10}^{11}\cdot \left( 1+\log l\right)   \label{3.11}
\end{equation}%
from (\ref{3.6}). Besides, using (\ref{1.2}), (\ref{3.0}) and (\ref{3.1}),
we can write\
\begin{equation}
0<\frac{(2^{k}-1)\cdot 3^{l}}{7^{m}}-1=\frac{2^{k}-2}{7^{m}}<\frac{2^{k}}{%
7^{m}}<\frac{4/3}{3^{l}}.  \label{3.12}
\end{equation}%
Put\noindent\ $\Lambda _{3}:=1-{{\gamma }_{1}}^{b_{1}}{{\gamma }_{2}}^{b_{2}}%
{{\gamma }_{3}}^{b_{3}}$, where%
\begin{equation*}
{\gamma }_{1}:=(2^{k}-1),{\gamma }_{2}:=3,{\gamma }%
_{3}:=7,b_{1}:=1,b_{2}:=l,b_{3}:=-m.
\end{equation*}%
A similar arguments to the above shows that $\Lambda _{3}\neq 0,$ Also, we
can take $D=1,$ $B:=l$, $A_{1}:=k{\log 2}$, $A_{2}:=\log 3$ and $A_{3}:={%
\log 7}$. Thus, applying Theorem \ref{T2} to inequality (\ref{3.12}), we
get\noindent
\begin{equation*}
l\log 3-\log (4/3)<1.4{\cdot 30}^{6}{\cdot }3^{4.5}\cdot \left( 1+\log
l\right) \cdot k\cdot {\log 2\cdot }\log 3\cdot {\log 7},
\end{equation*}%
which implies that
\begin{equation}
l<1.94\cdot 10^{11}\cdot \left( 1+\log l\right) \cdot k.  \label{3.13}
\end{equation}%
\noindent If we substitute the bound of $k$ in (\ref{3.11}) into (\ref{3.13}%
), then we get%
\begin{equation}
k<l<2.29\cdot 10^{26}.  \label{3.14}
\end{equation}%
Thus, we have the following lemma.

\begin{lemma}
\label{L7}If the equation (\ref{1.2}) holds with $k\geq 3$, $l\geq 2$ and $%
m\geq 3,$ then
\begin{equation}
k,l,m<3.72\cdot 10^{26}.  \label{3.15}
\end{equation}
\end{lemma}

Now, let reduce this upper bound on $k$ and $l.$ For this, let\noindent\
consider the linear form%
\begin{equation*}
z_{1}:=k\log 2-m\log 7+l\log {3.}
\end{equation*}%
Since this is a linear form in three integer variables, we use $LLL$%
-algorithm given in Lemma \ref{L4}. By Lemma \ref{L4}, say $X:=\max
\{X_{1},X_{2},X_{3}\},$ $S:=X_{1}^{2}+X_{2}^{2},~T:=\frac{1}{2}\left(
1+X_{1}+X_{2}+X_{3}\right) $ and $C>27X^{3}$ such that $k\leq X_{1},~l\leq
X_{2},$ and $m\leq X_{3}.$ In this case, by Lemma \ref{L7}, we can take
\begin{eqnarray*}
X &=&3.72\cdot 10^{26}, \\
S &=&3.76768\cdot 10^{53} \\
T &=&5.58\cdot 10^{26}
\end{eqnarray*}%
and thus $C=10^{82}.$ Moreover, the first vector of $LLL$-reduced basis is

\begin{equation*}
\left(
\begin{array}{c}
431878241430642097854554415 \\
-1794123608744301841719847442 \\
350763145311579067760845526%
\end{array}%
\right) .
\end{equation*}%
Further, we can see that $\left\Vert b_{1}\right\Vert \approx 1.87841\cdot
10^{27}$, $c_{1}=1,$ and thus $\theta =1.87841\cdot 10^{27}$ in Lemma \ref%
{L4}. Since $\theta ^{2}>S+T^{2},$ it follows that
\begin{equation}
\left\vert z_{1}\right\vert \geq \frac{\sqrt{\theta ^{2}-S}-T}{C}>1.245\cdot
10^{-55}  \label{3.16}
\end{equation}%
by Lemma \ref{L4}. Besides, if we say $\Lambda _{1}=1-e^{z_{1}},$ then, from
(\ref{3.5}), we get\noindent
\begin{equation*}
\left\vert \Lambda _{1}\right\vert =e^{z_{1}}-1<\frac{8/3}{2^{\min \left\{
k,l\right\} }}<\frac{8}{{12}}<0.7
\end{equation*}%
since $k\geq 3$ and $l\geq 2.$ Taking $a=0.7$ in Lemma \ref{L3}, it is seen
that

\noindent
\begin{equation*}
|z_{1}|=|\log (\Lambda _{1}+1)|<\frac{{-\log \left( 1-0.7\right) \ }}{0.7}%
\cdot \left\vert \Lambda _{1}\right\vert <\frac{5}{2^{\min \left\{
k,l\right\} }}.
\end{equation*}%
Hence, from (\ref{3.16}), we have%
\begin{equation*}
1.245\cdot 10^{-55}<\frac{5}{2^{\min \left\{ k,l\right\} }},
\end{equation*}%
which implies that $\min \left\{ k,l\right\} \leq 184.$

Assume that $l\leq k.$ Then $l\leq 184.$ Substituting this bound of $l$ into
(\ref{3.9}), we get%
\begin{equation*}
k<3.1\cdot 10^{11}\cdot \left( 1+\log k\right) {\cdot 184},
\end{equation*}%
which implies that $k<2.069\cdot 10^{15}$. Now, to further reduce this bound
on $k,$ we apply Lemma \ref{L2}. Then, let\noindent
\begin{equation*}
z_{2}:=k\log 2-m\log 7+\log {\left( 3^{l}-1\right) }
\end{equation*}%
and $\Lambda _{2}=1-e^{z_{2}}.$ Thus, from (\ref{3.8}), it is seen
that\noindent
\begin{equation*}
\left\vert \Lambda _{2}\right\vert =e^{z_{2}}-1<\frac{4/3}{2^{k}}<\frac{1}{6}%
<0.2
\end{equation*}%
since $k\geq 3.$ Choosing $a=0.2$ in Lemma \ref{L3}, we obtain\noindent
\begin{equation*}
|z_{2}|=|\log (\Lambda _{2}+1)|<\frac{{-\log \left( 1-0.2\right) \ }}{0.2}%
\cdot \left\vert \Lambda _{2}\right\vert <\frac{1.5}{2^{k}},
\end{equation*}%
that is,%
\begin{equation*}
\left\vert k\log 2-m\log 7+\log {\left( 3^{l}-1\right) }\right\vert <\frac{%
1.5}{2^{k}}.
\end{equation*}%
Dividing this inequality by $\log 7$, we get

\noindent
\begin{equation}
\left\vert k\gamma -m+\mu \ \right\vert <A\cdot B^{-w},  \label{3.17}
\end{equation}%
where\noindent \noindent
\begin{equation*}
\gamma =\frac{\log 2}{\log 7}\text{, }\mu =\frac{\log {\left( 3^{l}-1\right)
\ }}{\log 7}\text{, }A=0.78\text{,}\ B=2\text{ and }w=k.
\end{equation*}%
$\gamma $ is a irrational number. Because, if it were not, it would be $%
\gamma =\frac{\log 2}{\log 7}=\frac{a}{b}$ for some integers $a$ and $b.$
This implies that $2^{b}=7^{b},$ which is impossible. Taking $M:=2.069\cdot
10^{15},$ which is an upper bound on $k$, we find that the denominator of
the $30$ th convergent of $\gamma $ exceeds $6M.$ Also, a quick computation
with Mathematica gives that the inequality%
\begin{equation*}
0<0.000190128<\varepsilon (l)
\end{equation*}%
is valid for all values of $2\leq l\leq 184,$ where $l$ is odd. Thus, a
similar argument to the above case shows that if inequality (\ref{3.17}) has
a solution, then%
\begin{equation*}
w=k<68.555.
\end{equation*}%
Thus, since $l\leq k,$ it follows that $l\leq 68.$

Now, assume that $k<l.$ Then $k\leq 184.$ Substituting this bound of $k$
into (\ref{3.13}), we get%
\begin{equation*}
l<1.94\cdot 10^{11}\cdot \left( 1+\log l\right) \cdot 184,
\end{equation*}%
which implies that $l<1.278\cdot 10^{15}$. Now, to further reduce this bound
on $l,$ we apply Lemma \ref{L2}. Then, let\noindent
\begin{equation*}
z_{3}:=l\log 3-m\log 7+\log {\left( 2^{k}-1\right) }
\end{equation*}%
and $\Lambda _{3}=1-e^{z_{3}}.$ Thus, from (\ref{3.12}), it is seen
that\noindent
\begin{equation*}
\left\vert \Lambda _{3}\right\vert =e^{z_{3}}-1<\frac{4/3}{3^{l}}\leq \frac{4%
}{27}<0.2
\end{equation*}%
since $l\geq 2.$ Choosing $a=0.2$ in Lemma \ref{L3}, we obtain\noindent
\begin{equation*}
|z_{3}|=|\log (\Lambda _{3}+1)|<\frac{{-\log \left( 1-0.2\right) \ }}{0.2}%
\cdot \left\vert \Lambda _{3}\right\vert <\frac{1.5}{3^{l}},
\end{equation*}%
that is,%
\begin{equation*}
\left\vert l\log 3-m\log 7+\log {\left( 2^{k}-1\right) }\right\vert <\frac{%
1.5}{3^{l}}.
\end{equation*}%
Dividing this inequality by $\log 7$, we get

\noindent
\begin{equation}
\left\vert l\gamma -m+\mu \ \right\vert <A\cdot B^{-w},  \label{3.18}
\end{equation}%
where\noindent \noindent
\begin{equation*}
\gamma =\frac{\log 3}{\log 7}\text{, }\mu =\frac{\log {\left( 2^{k}-1\right)
\ }}{\log 7}\text{, }A=0.78\text{,}\ B=3\text{ and }w=l.
\end{equation*}%
It can be seen similarly to the above that $\gamma =\frac{\log 3}{\log 7}$
is a irrational number. Taking $M:=1.278\cdot 10^{15},$ which is an upper
bound on $l$, we find that the denominator of the $34$ th convergent of $%
\gamma $ exceeds $6M.$ Also, a quick computation with Mathematica gives that
$\varepsilon (k)$ is positive for all $3\leq k\leq 184$ except for $%
\varepsilon (3)$ and $\varepsilon (6)$. Let $k\neq 3$ and $k\neq 6.$ In this
case, Lemma \ref{L2} says that if inequality (\ref{3.18}) has a solution,
then%
\begin{equation*}
w=l<44.041.
\end{equation*}%
Further, if $k=3$ or $k=6,$ then we have $7\cdot {\left( 3^{l}-1\right) =7}%
^{m}-1$ or $63\cdot {\left( 3^{l}-1\right) =7}^{m}-1.$ But these are
impossible. Thus, since $k<l,$ it follows that $k\leq 44.$

Consequently, we can give the following lemma.

\begin{lemma}
\label{L8}The equation (\ref{1.2}) holds with $k\geq 3$, $l\geq 2,$ and $%
m\geq 3.$ Then, $\max \left\{ k,l\right\} \leq 68.$
\end{lemma}

Lastly, a simple computation with Mathematica shows that equation (\ref{3.18}%
) has no solutions for $3\leq k\leq 68,~2\leq l\leq 68$ and $3\leq m<68.$
Thus, the proof of Theorem \ref{T1} is completed.

In our next study, we will investigate the positive integer quadruples $%
(n,k,l,m)$ satisfying the equation
\begin{equation}
(L_{n}^{k}-1)(L_{n+1}^{l}-1)=L_{n+2}^{m}-1  \label{3.19}
\end{equation}%
that generalizes the equation (\ref{1.2}). In (\ref{3.19}), taking $n=2,$ we
have $(3^{k}-1)(4^{l}-1)=7^{m}-1.$ This equation has only the solution $%
(k,l,m)=(1,1,1)$ by Theorem \ref{T1}.

\end{document}